\documentclass[11pt]{article}

\usepackage[T1]{fontenc}
\usepackage[utf8]{inputenc}
\usepackage[margin=1in]{geometry}
\usepackage{amsmath,amssymb,amsfonts,amsthm,bm}
\usepackage{graphicx}
\usepackage{booktabs}
\usepackage{xcolor}
\usepackage{hyperref}
\usepackage{mathtools}
\usepackage{enumitem}
\usepackage{subcaption}
\usepackage{float}
\usepackage{authblk}
\usepackage{lineno}
\usepackage{microtype}

\graphicspath{{figures/}}

\hypersetup{
    colorlinks=true,
    linkcolor=blue,
    citecolor=blue,
    urlcolor=blue
}

\title{
Structural Causal Discovery and Predictive Sufficiency in High-Dimensional Dynamical Systems
}

\author[1]{Abd AlRahman R. AlMomani}
\author[1]{Curtis N. James}
\author[2]{Christopher C. Hennon}
\author[1]{Ronny Schroeder}

\affil[1]{Embry-Riddle Aeronautical University, Prescott, Arizona, 86301}
\affil[2]{Western Carolina University, Cullowhee, North Carolina, 28723}
\date{}

\begin{document}

\maketitle

\section*{Highlights}
\label{front:highlights}

\begin{itemize}
    \item Entropic regression identifies a compact, spatially recurrent set of candidate causal parents of next-hour precipitation.
    \item A variable-level forward and backward formulation discovers nonlinear structure without constructing the global quadratic library.
    \item Under a common cardinality, entropic regression yields a more concentrated causal-parent recurrence profile than transfer entropy and causation entropy.
    \item The recovered causal structure retains strong one-hour occurrence discrimination, while intensity prediction remains more limited.
\end{itemize}

\section*{eTOC Blurb}
\label{front:etoc}

AlMomani et al. develop a projection-based, variable-level formulation of entropic regression for structural causal discovery in high-dimensional precipitation dynamics. Across 1227 spatial regions, the method recovers a compact and physically coherent set of recurrent candidate causal parents, with a more concentrated recurrence profile than transfer entropy and causation entropy. Predictive tests show that this causal structure retains strong short-horizon occurrence information while resolving precipitation intensity less completely.

\begin{abstract}
High-dimensional environmental systems contain variables that may be predictive, causally informative, physically coupled, or redundant, and these roles need not coincide. We examine this distinction in precipitation dynamics using High-Resolution Rapid Refresh atmospheric fields and Multi-Radar Multi-Sensor precipitation observations over the Southwestern United States. The principal methodological contribution is a projection-based, variable-level formulation of entropic regression that identifies candidate causal parents without constructing a global nonlinear expansion library. Ten physical variables are admitted sequentially through conditional-information forward selection using linear regression projections. A quadratic library is then generated only over the selected variables, and two sequential variable-block eliminations produce the final eight-variable local representation. Applied across 1227 precipitation-guided superpixels, the method recovers a compact and spatially recurrent causal structure. Under the same final cardinality, entropic regression produces a substantially more concentrated recurrence profile than transfer entropy and causation entropy. Six variables exceed an operational recurrence threshold of 0.25: land-surface moisture availability, composite reflectivity, geometric vertical velocity, wind speed, convective available potential energy, and maximum upward vertical velocity. Predictive models restricted to these variables achieve strong one-hour discrimination of precipitation occurrence, with an area under the receiver operating characteristic curve of 0.948, while pointwise intensity prediction remains substantially more limited. These results establish two contributions: a scalable variable-level entropic regression procedure for high-dimensional causal discovery, and empirical evidence that recurrent causal relevance and predictive sufficiency are distinct properties of a reduced dynamical representation. Within the explicit observational definition adopted in this paper, the six recurrent variables constitute the causal backbone recovered by entropic regression.
\end{abstract}

\section{Introduction}
\label{sec:introduction}

Understanding the organization of precipitation remains a central challenge in atmospheric and Earth-system science. Rainfall emerges through coupled thermodynamic, kinematic, microphysical, and land-surface processes, with strong dependence on moisture, instability, lifting, cloud microphysics, flow organization, and terrain \cite{Houze2014,Markowski2010,Koster2004,Seneviratne2010}. These interactions are especially heterogeneous during the North American monsoon, whose precipitation varies strongly across elevation, regional circulation, and convective regime \cite{Adams1997Monsoon}. The resulting observations are sparse, intermittent, nonlinear, and spatially dependent. A variable may be associated with precipitation because it participates in the evolving process, because it is a proxy for another field, because it carries persistence information, or because it shares an unresolved common driver with the response.

Numerical weather prediction provides a physically organized framework for atmospheric forecasting \cite{Kalnay2003,dowell2022hrrr}. In parallel, data-driven methods have become increasingly important for weather and Earth-system applications \cite{reichstein2019deep,Schultz2021,rasp2020weatherbench,rasp2024weatherbench2}. Convolutional, recurrent, and generative models can approximate complex predictive mappings and have produced important advances in precipitation nowcasting \cite{ayzel2020rainnet,shi2015convlstm,ravuri2021skilful}. Predictive performance, however, does not by itself determine which inputs carry nonredundant structural information. A flexible model may exploit a large collection of correlated fields without isolating a compact representation of the process.

The reverse implication is also not guaranteed. A variable set selected because it contributes unique information under a structural criterion may not form a sufficient predictive state for every forecast objective. Forecast quality is multidimensional: discrimination, reliability, sharpness, threshold-dependent event detection, quantitative intensity, and spatial organization answer different questions \cite{Murphy1993,Gneiting2007,Fawcett2006,Saito2015,Gilleland2009,Wilks2019}. A reduced representation may perform well for one objective while remaining incomplete for another. This distinction is central to the present paper.

For observational time series, causal discovery must be defined through an explicit information structure and a stated set of assumptions. Latent variables, indirect pathways, measurement error, contemporaneous dependence, and incomplete state description remain central concerns \cite{Pearl2009,Peters2017,Runge2019CausationEarth}. In this study, lagged incremental conditional information provides the operational criterion for identifying candidate causal parents. This criterion is stronger than marginal association because every admitted variable must contribute information beyond the variables already selected, while remaining distinct from intervention-level identification.

Information-theoretic methods provide a natural language for separating shared and incremental dependence. Transfer entropy measures directed information beyond the past of the target \cite{Schreiber2000,Barnett2009}, while causation entropy extends conditional-information screening to multivariate parent-set discovery \cite{Sun2015OCSE}. Their practical application in high-dimensional continuous systems remains difficult because nearest-neighbor estimates of mutual information and conditional mutual information are sensitive to sample size, dimension, temporal dependence, repeated values, and mixed discrete-continuous structure \cite{kraskov2004estimating,frenzel2007partial,Runge2018CMI,Gao2017Mixture,Mesner2021MixedCMI}.

Entropic regression was introduced as an information-theoretic approach to nonlinear system identification \cite{almomani2020entropic,fish2021entropic,almomani2020erfit}. It combines iterative model construction with conditional-information selection. In the standard formulation, one basis function is admitted at each forward iteration according to the information carried by its regression projection after conditioning on the current model. A backward stage then reexamines the admitted terms and removes those that are redundant. This logic is related to optimal causation entropy, but the evaluated objects are regression projections rather than the raw candidate variables. It also differs from sparse-library methods such as SINDy, which begin from an explicitly constructed candidate function library \cite{Brunton2016SINDy}.

The present work adapts this principle to a predictor space containing 122 physical variables. The forward stage operates only on the linear variable terms and selects ten variables sequentially. A quadratic expansion is then generated over this reduced set, after which two sequential backward eliminations remove variables together with their associated nonlinear terms. Thus, the full quadratic library over all 122 variables is never constructed. This modification preserves the forward and backward organization of entropic regression while keeping the selection unit physically interpretable.

The spatial analysis is also designed as a representation problem. Data-driven partitions and spatially coherent structures have previously been used to organize image-observed dynamics, remote-sensing fields, and spatiotemporal transitions \cite{almomani2020go,almomani2021early}. Related work on geometric partition entropy and its generalizations emphasizes that finite partitions and geometry materially affect information extracted from continuous data \cite{diggans2022geometric,diggans2023boltzmann,diggans2025generalizing}. Here, precipitation-guided superpixels define a fixed collection of local causal-discovery problems, and the recurrence of selected parents across these regions provides the domain-level summary.

\paragraph{Contributions.}
This paper makes two principal contributions. Methodologically, it introduces a projection-based, variable-level formulation of entropic regression that identifies causal structure in a high-dimensional nonlinear system without constructing a global expansion library. Empirically, it recovers a compact and physically coherent causal backbone whose conditional information contribution persists across the precipitation domain. Predictive experiments then show that this backbone retains strong next-hour occurrence information while resolving precipitation intensity less completely. The resulting separation between causal relevance and predictive sufficiency is a positive structural finding, not merely a limitation of the predictive model.

The analysis has four objectives. First, we formulate the variable-level entropic regression procedure for lagged structural causal discovery. Second, we apply the procedure independently over precipitation-guided spatial regions and summarize the recurrence of selected causal parents across the domain. Third, we compare the resulting recurrence profile with those produced by transfer entropy and causation entropy under a common final cardinality. Fourth, we use the recurrent causal structure in predictive models to characterize how much occurrence, intensity, and spatial information is retained by the reduced representation.

Throughout this study, the phrase \emph{candidate causal parent} denotes a lagged variable that contributes incremental information about next-hour precipitation after conditioning on the currently selected model, under the adopted temporal pairing, regression operator, information estimator, cardinality, and spatial representation. After this operational definition is established, the terms \emph{causal structure} and \emph{causal backbone} refer to the recurrent parent set recovered by that procedure. Predictive sufficiency is evaluated relative to the variables, history windows, model architectures, training procedures, and metrics considered here.

The remainder of the paper is organized as follows. Section~\ref{sec:problem_formulation} states the structural and predictive questions. Section~\ref{sec:data_spatial} describes the data and spatial representation. Section~\ref{sec:entropic_regression} develops the modified entropic regression procedure. Section~\ref{sec:comparison} introduces the transfer entropy and causation entropy comparison. Section~\ref{sec:results} reports the structural and predictive results. Section~\ref{sec:discussion} discusses their interpretation and limitations.

\section{Problem Formulation}
\label{sec:problem_formulation}

Let $Y(\mathbf{s},t)$ denote precipitation at spatial location $\mathbf{s}$ and time $t$, and let
\begin{equation}
\mathbf{X}(\mathbf{s},t)
=
\big(X_1(\mathbf{s},t),\ldots,X_p(\mathbf{s},t)\big)
\label{eq:predictor_vector}
\end{equation}
denote the available atmospheric predictors. The structural analysis is lagged by one hour: predictors observed at time $t$ are evaluated with respect to precipitation at time $t+1$.

The study addresses two related questions:
\begin{enumerate}[label=(\roman*)]
    \item Which variables act as recurrent candidate causal parents of $Y(\mathbf{s},t+1)$ by providing nonredundant lagged information across the adopted spatial partition?
    \item What predictive information is retained when the recovered causal structure is used as the only input to occurrence and intensity models?
\end{enumerate}

The first question is a structural causal-discovery problem. It asks whether sequential conditional-information selection resolves a compact parent structure that persists across heterogeneous local representations. The second is a predictive characterization problem. It asks which forecast attributes are retained by that causal structure, rather than assuming that parent discovery and optimized prediction are equivalent tasks.

\section{Data and Spatial Representation}
\label{sec:data_spatial}

\subsection{Data sources}

Let $\mathbf{s}\in\Omega\subset\mathbb{R}^2$ denote spatial location and let $t$ index the hourly observations retained for the North American monsoon seasons considered in the study. The response field is
\begin{equation}
Y(\mathbf{s},t),
\label{eq:response_field}
\end{equation}
representing hourly precipitation accumulation, and the predictor vector is
\begin{equation}
\mathbf{X}(\mathbf{s},t)
=
\big(X_1(\mathbf{s},t),\ldots,X_p(\mathbf{s},t)\big),
\qquad p=122.
\label{eq:predictor_field}
\end{equation}

The predictor fields are obtained from HRRR analysis data \cite{dowell2022hrrr,NOAAHRRRInventory}, and the precipitation response is obtained from MRMS observations \cite{zhang2016mrms}. The retained timestamps sample the North American monsoon season, a regime characterized by strong spatial and temporal variability across the Southwestern United States \cite{Adams1997Monsoon}. For each year $y\in\{2018,2019,2020,2021,2022\}$, let $\mathcal{T}_y$ denote the set of hourly timestamps included from that year's monsoon-season sample. The development and deployment periods are then
\begin{equation}
\mathcal{T}_{\mathrm{dev}}
=
\mathcal{T}_{2018}\cup\mathcal{T}_{2019}\cup\mathcal{T}_{2020}\cup\mathcal{T}_{2021},
\label{eq:development_period}
\end{equation}
and
\begin{equation}
\mathcal{T}_{\mathrm{test}}=\mathcal{T}_{2022}.
\label{eq:test_period}
\end{equation}
All structural screening and model development use $\mathcal{T}_{\mathrm{dev}}$. The 2022 sample is reserved for deployment-level evaluation.

HRRR is an hourly updating, convection-allowing system with approximately $3\,\mathrm{km}$ horizontal grid spacing and radar data assimilation in its operational analysis cycle \cite{dowell2022hrrr}. The retained fields have hourly temporal resolution,
\begin{equation}
\Delta t=1\,\mathrm{hour}.
\label{eq:hrrr_resolution}
\end{equation}
MRMS observations are aggregated to hourly accumulation and resampled to the HRRR grid for predictive modeling.

Two representations are used for different purposes:
\begin{itemize}
    \item \textbf{Structural screening representation:} predictors and precipitation are averaged within spatially contiguous superpixels.
    \item \textbf{Predictive representation:} MRMS precipitation and HRRR predictors are retained as aligned grid fields.
\end{itemize}
The superpixel representation is used only for structural screening and domain-level recurrence analysis. Predictive models are trained and evaluated on the aligned grid representation.

\subsection{High-dimensional structure}
\label{sec:high_dimensional_structure}

At each time $t$, the discretized system is represented by
\begin{equation}
\mathbf{X}(t)\in\mathbb{R}^{N\times p},
\qquad
Y(t)\in\mathbb{R}^{N},
\label{eq:grid_representation}
\end{equation}
where $N$ is the number of grid locations. Three properties motivate the structural representation used below.

\paragraph{Spatial redundancy.}
Neighboring locations are strongly dependent. Treating every grid location as an independent inference unit would therefore overstate the effective spatial sample size and repeat highly similar local problems.

\paragraph{Variable coupling.}
Many HRRR fields describe related atmospheric processes or different representations of the same evolving state. Their explanatory contributions may overlap substantially, so marginal association alone is not sufficient for identifying a compact set.

\paragraph{Precipitation sparsity.}
Hourly precipitation is zero over much of the space-time domain. Define
\begin{equation}
\mathcal{I}_0
=
\{(\mathbf{s},t):Y(\mathbf{s},t)=0\},
\qquad
\mathcal{I}_1
=
\{(\mathbf{s},t):Y(\mathbf{s},t)>0\}.
\label{eq:index_sets}
\end{equation}
Then
\begin{equation}
|\mathcal{I}_0|\gg|\mathcal{I}_1|.
\label{eq:imbalance}
\end{equation}
The zero regime is scientifically meaningful, but event-specific information is concentrated in the much smaller nonzero subset.

\begin{figure}[t]
\centering
\includegraphics[width=0.78\columnwidth]{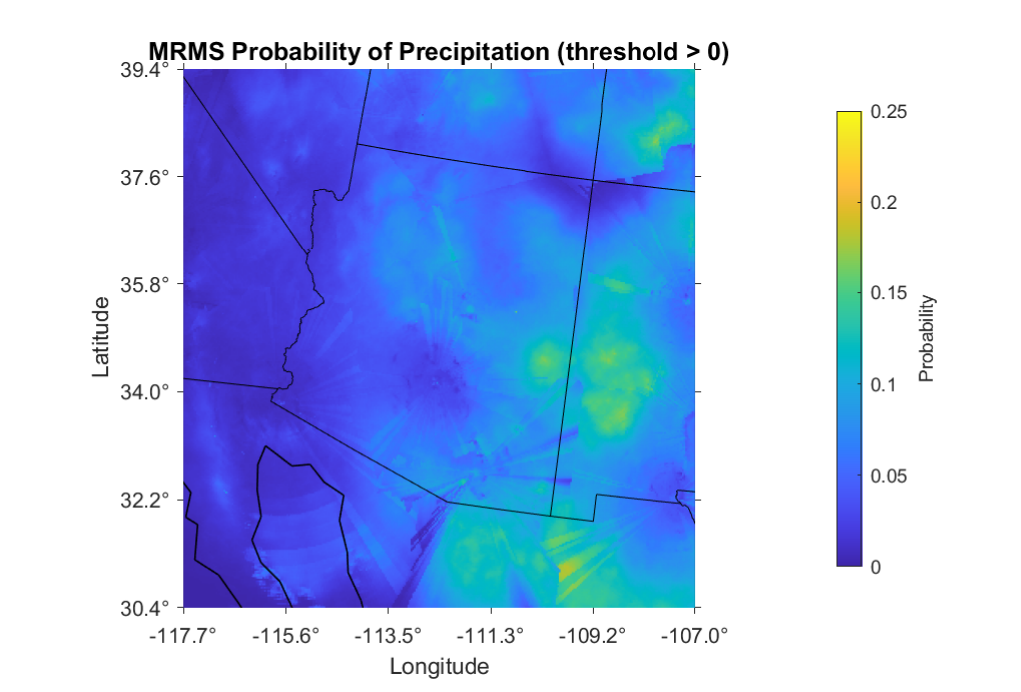}
\caption{Spatial distribution of the empirical precipitation probability $P(\mathbf{s})$ over the development period. The field shows broad regions of low occurrence together with localized regions of higher monsoon-season activity.}
\label{fig:pop_map}
\end{figure}

The empirical probability of precipitation at location $\mathbf{s}$ is
\begin{equation}
P(\mathbf{s})
=
\frac{1}{|\mathcal{T}_{\mathrm{dev}}|}
\sum_{t\in\mathcal{T}_{\mathrm{dev}}}
\mathbf{1}_{\{Y(\mathbf{s},t)>0\}},
\label{eq:pop}
\end{equation}
where $\mathbf{1}_{\{\cdot\}}$ is the indicator function. Figure~\ref{fig:pop_map} shows that precipitation occurrence is spatially heterogeneous over the study domain.

For later interpretation, define the number of nonzero precipitation observations at location $\mathbf{s}$ by
\begin{equation}
N_{\mathrm{rain}}(\mathbf{s})
=
\sum_{t\in\mathcal{T}_{\mathrm{dev}}}
\mathbf{1}_{\{Y(\mathbf{s},t)>0\}}.
\label{eq:neff}
\end{equation}
The quantity $N_{\mathrm{rain}}(\mathbf{s})$ is not an effective sample size in the formal dependence-adjusted sense. It is a descriptive count of event observations and emphasizes that the amount of rain-regime information varies across the domain.

\subsection{Superpixel representation}
\label{sec:superpixel_representation}

For each timestamp $t_k$, the precipitation field is vectorized as
\begin{equation}
\mathbf{y}_k\in\mathbb{R}^{N},
\qquad
(\mathbf{y}_k)_i=Y(\mathbf{s}_i,t_k).
\label{eq:vectorization}
\end{equation}
Stacking these vectors gives
\begin{equation}
\mathbf{Y}
=
\begin{bmatrix}
\mathbf{y}_1&\mathbf{y}_2&\cdots&\mathbf{y}_T
\end{bmatrix}
\in\mathbb{R}^{N\times T}.
\label{eq:data_matrix}
\end{equation}
The temporal mean is removed at each spatial location,
\begin{equation}
\widetilde{\mathbf{Y}}
=
\mathbf{Y}-\overline{\mathbf{Y}},
\qquad
\overline{\mathbf{Y}}_{i,k}
=
\frac{1}{T}\sum_{k'=1}^{T}Y(\mathbf{s}_i,t_{k'}),
\label{eq:centering}
\end{equation}
and the singular value decomposition, used here as a descriptive low-rank representation, is \cite{Jolliffe2016PCA}
\begin{equation}
\widetilde{\mathbf{Y}}
=
\mathbf{U}\boldsymbol{\Sigma}\mathbf{V}^{\top}.
\label{eq:svd}
\end{equation}
The leading spatial mode is represented by
\begin{equation}
\mathrm{PC}_1=\sigma_1\mathbf{u}_1\in\mathbb{R}^{N}.
\label{eq:pc1}
\end{equation}
We define the scalar partitioning field
\begin{equation}
\phi(\mathbf{s}_i)=(\mathrm{PC}_1)_i.
\label{eq:phi}
\end{equation}

Simple Linear Iterative Clustering is then applied to $\phi$ to obtain spatially contiguous regions that balance feature similarity and spatial proximity \cite{achanta2012slic}. This yields the fixed partition
\begin{equation}
\Omega
=
\bigcup_{m=1}^{M}\Omega_m,
\qquad
\Omega_m\cap\Omega_{m'}=\emptyset
\quad(m\neq m'),
\qquad
M=1227.
\label{eq:partition}
\end{equation}

For each superpixel, the aggregated response and predictor series are
\begin{equation}
\overline{Y}_m(t)
=
\frac{1}{|\Omega_m|}
\sum_{\mathbf{s}_i\in\Omega_m}Y(\mathbf{s}_i,t),
\label{eq:agg_y}
\end{equation}
and
\begin{equation}
\overline{X}_{m,j}(t)
=
\frac{1}{|\Omega_m|}
\sum_{\mathbf{s}_i\in\Omega_m}X_j(\mathbf{s}_i,t),
\qquad j=1,\ldots,p.
\label{eq:agg_x}
\end{equation}
The resulting multivariate series $(\overline{\mathbf{X}}_m(t),\overline{Y}_m(t))$ defines one local structural screening problem. Spatial coarse-graining is used here as a representation strategy rather than as a claim that the superpixels are independent replicates. Related work has used data-driven spatial partitions to extract coherent structure from image-observed dynamics and remote-sensing fields \cite{almomani2020go,almomani2021early}. Geometric partition entropy, Boltzmann-Shannon interaction entropy, and generalized geometric partition entropy further emphasize that partition geometry and sample support affect information extracted from continuous data \cite{diggans2022geometric,diggans2023boltzmann,diggans2025generalizing}.

\begin{figure}[t]
\centering
\includegraphics[width=0.96\columnwidth]{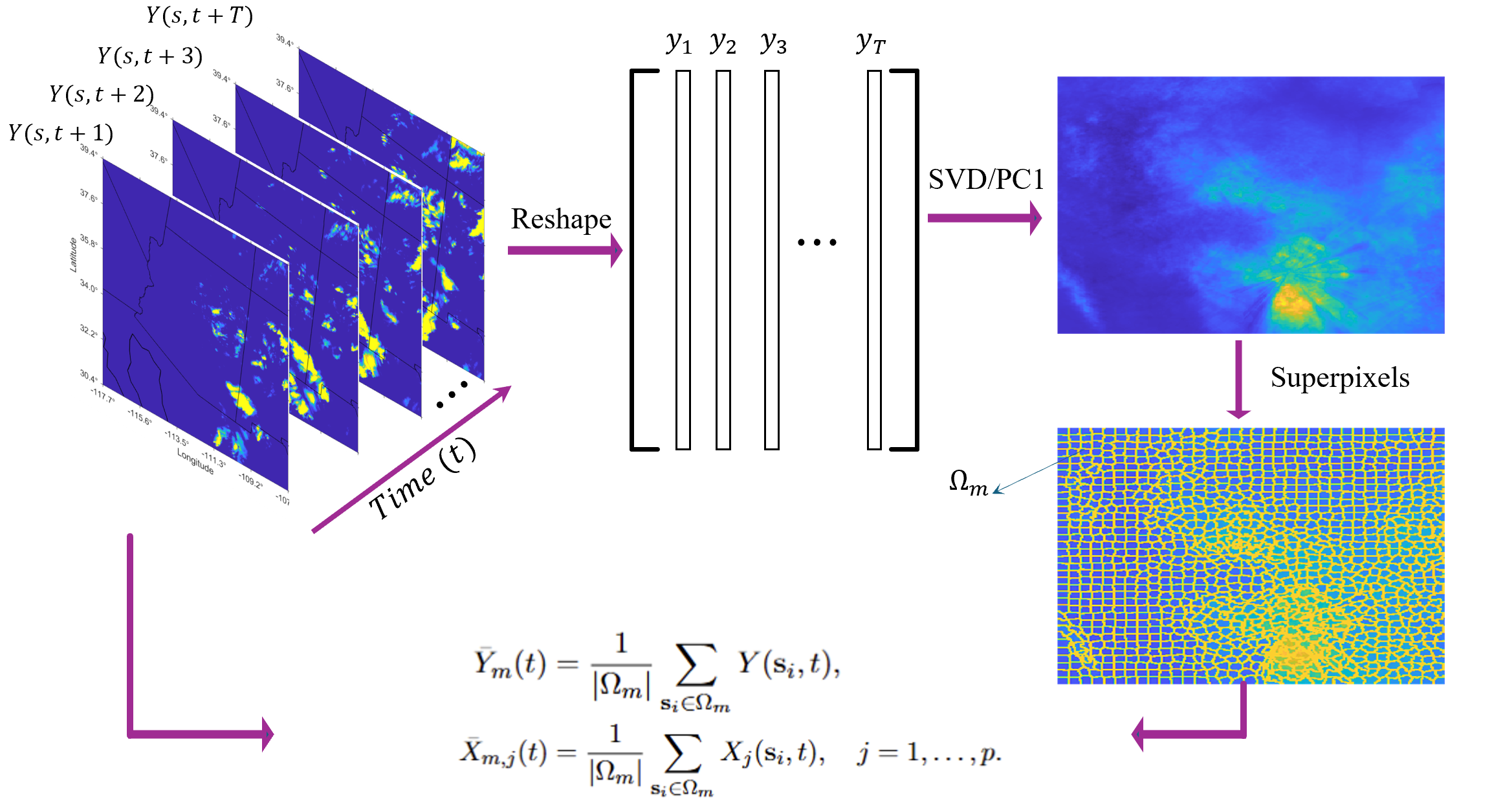}
\caption{Superpixel representation for structural screening. A dominant spatial mode of the precipitation sequence guides the SLIC partition. Atmospheric predictors and precipitation are averaged within each region to obtain the local time series used by the information-theoretic selection procedures.}
\label{fig:superpixel_representation}
\end{figure}

This construction emphasizes the dominant spatial organization of the development data. Localized or transient structures that are weakly represented in the first mode may be attenuated. The consequence is addressed as a limitation rather than hidden within the recurrence statistic.

\section{Entropic Regression for Variable-Level Structural Causal Discovery}
\label{sec:entropic_regression}

\subsection{Motivation and relation to standard ER}

Entropic regression combines regression projections with iterative information-theoretic selection \cite{almomani2020entropic,fish2021entropic,almomani2020erfit}. In the standard formulation, a candidate library is specified before selection. At each forward iteration, one library term is admitted according to the conditional information carried by its regression projection after conditioning on the projection from the currently selected terms. The backward stage reexamines the forward set and removes the term with the smallest conditional contribution. The same basic organization underlies the variable-level modification used here.

For $p$ physical predictors, a quadratic library containing all linear terms, squares, and unique pairwise products has
\begin{equation}
D(p)=p+\frac{p(p+1)}{2}
\label{eq:global_library_dimension}
\end{equation}
columns. For $p=122$, this gives $D(122)=7625$ candidate functions before any screening. Large preconstructed libraries are common in sparse nonlinear identification \cite{Brunton2016SINDy}, but they mix physical-variable identity with interaction identity and can be costly when a full selection procedure must be repeated across 1227 spatial regions.

The present modification separates variable screening from nonlinear expansion. The forward stage is performed on the 122 linear physical variables. After ten variables have been selected sequentially, a quadratic library is constructed only over those variables. The backward stage then removes variables at the group level, meaning that removal of a variable also removes its square and every interaction involving that variable. This design has three defining properties:
\begin{enumerate}[label=(\roman*)]
    \item physical variables, rather than individual polynomial terms, are the selection units;
    \item the forward stage preserves the sequential conditional-information structure of ER but uses linear variable projections;
    \item nonlinear terms are generated only after forward screening and are refined through variable-block backward elimination.
\end{enumerate}
Together, these changes turn the original library-based ER logic into a variable-level causal-discovery procedure whose output remains interpretable in terms of physical atmospheric fields.

\subsection{Lagged local response and projection operator}

For superpixel $\Omega_m$, define the lagged target vector
\begin{equation}
\mathbf{y}_m
=
\begin{bmatrix}
\overline{Y}_m(2)&\overline{Y}_m(3)&\cdots&\overline{Y}_m(T)
\end{bmatrix}^{\top}
\in\mathbb{R}^{T-1}.
\label{eq:target_vector}
\end{equation}
For a set of variable indices $\mathcal{A}$, the linear design matrix is
\begin{equation}
\mathbf{L}_m(\mathcal{A})
=
\left[
\overline{X}_{m,j}(t)
\right]_{
 t=1,\ldots,T-1;\,j\in\mathcal{A}}
\in\mathbb{R}^{(T-1)\times|\mathcal{A}|}.
\label{eq:linear_design}
\end{equation}
For any design matrix $\mathbf{F}$ with $T-1$ rows, define the least-squares projection of the target by
\begin{equation}
\mathcal{P}_m(\mathbf{F})
=
\mathbf{F}\mathbf{F}^{\dagger}\mathbf{y}_m,
\label{eq:projection_operator}
\end{equation}
where $\mathbf{F}^{\dagger}$ denotes the Moore-Penrose pseudoinverse. This notation isolates the central ER object: the information criterion is evaluated on fitted response projections, not directly on the full raw predictor vector.

The local structural model is written conceptually as
\begin{equation}
\overline{Y}_m(t+1)
=
f_m\!\left(\{\overline{X}_{m,j}(t)\}_{j\in\mathcal{S}_m}\right)
+\eta_m(t+1),
\label{eq:nonlinear_model}
\end{equation}
where $\mathcal{S}_m$ is the final variable set for region $m$. The one-hour lag is fixed throughout the structural analysis.

\subsection{Sequential linear forward selection}
\label{sec:forward_selection}

The forward stage starts from
\begin{equation}
\mathcal{S}^{(0)}_m=\emptyset.
\label{eq:init_set}
\end{equation}
At the first iteration, each candidate variable $j$ is represented by the one-variable projection
\begin{equation}
\mathbf{z}_{m,j}
=
\mathcal{P}_m\!\left(\mathbf{L}_m(\{j\})\right),
\label{eq:candidate_linear_projection}
\end{equation}
and receives the score
\begin{equation}
\Delta I^{(1)}_{m,j}
=
I\!\left(\mathbf{y}_m;\mathbf{z}_{m,j}\right).
\label{eq:forward_first}
\end{equation}
The selected index is
\begin{equation}
j_1
=
\arg\max_{j\in\{1,\ldots,p\}}
\Delta I^{(1)}_{m,j},
\qquad
\mathcal{S}^{(1)}_m=\{j_1\}.
\label{eq:forward_first_update}
\end{equation}

For forward iteration $r\geq2$, the current selected-model projection is
\begin{equation}
\mathbf{z}_{m,\mathcal{S}^{(r-1)}_m}
=
\mathcal{P}_m\!\left(\mathbf{L}_m(\mathcal{S}^{(r-1)}_m)\right).
\label{eq:current_linear_projection}
\end{equation}
Every unselected candidate is then scored by
\begin{equation}
\Delta I^{(r)}_{m,j}
=
I\!\left(
\mathbf{y}_m;
\mathbf{z}_{m,j}
\,\middle|\,
\mathbf{z}_{m,\mathcal{S}^{(r-1)}_m}
\right),
\qquad
j\notin\mathcal{S}^{(r-1)}_m.
\label{eq:deltaI}
\end{equation}
The maximizing candidate is admitted,
\begin{equation}
j_r
=
\arg\max_{j\notin\mathcal{S}^{(r-1)}_m}
\Delta I^{(r)}_{m,j},
\qquad
\mathcal{S}^{(r)}_m
=
\mathcal{S}^{(r-1)}_m\cup\{j_r\}.
\label{eq:update_forward}
\end{equation}
In this study, $K_f=10$. Thus, the procedure performs ten sequential forward iterations. One variable is admitted at each iteration, and every score after the first is conditioned on the projection formed by the variables selected in all preceding iterations. The value $K_f$ is not the number of variables admitted simultaneously.

\subsection{Quadratic expansion after forward screening}

Let
\begin{equation}
\mathcal{B}^{(0)}_m
=
\mathcal{S}^{(K_f)}_m
\label{eq:forward_final_set}
\end{equation}
be the ten-variable forward set. Only at this stage is the quadratic design matrix generated:
\begin{equation}
\mathbf{Q}_m(\mathcal{B}^{(0)}_m)
=
\left[
\overline{X}_{m,i}(t),
\overline{X}_{m,i}(t)\overline{X}_{m,j}(t)
\;\middle|\;
 i,j\in\mathcal{B}^{(0)}_m,\ i\leq j
\right]_{t=1}^{T-1}.
\label{eq:quadratic_design}
\end{equation}
For a variable set $\mathcal{B}$, the number of columns is
\begin{equation}
d(\mathcal{B})
=
|\mathcal{B}|+\frac{|\mathcal{B}|(|\mathcal{B}|+1)}{2}.
\label{eq:feature_dimension}
\end{equation}
Consequently, the ten-variable forward set produces 65 linear and quadratic terms, rather than the 7625 terms in the global quadratic library.

\subsection{Sequential variable-block backward elimination}
\label{sec:backward_elimination}

At backward iteration $q$, let $\mathcal{B}^{(q-1)}_m$ denote the current variable set. For $j\in\mathcal{B}^{(q-1)}_m$, define the reduced set
\begin{equation}
\mathcal{B}^{(q-1)}_{m,-j}
=
\mathcal{B}^{(q-1)}_m\setminus\{j\}.
\label{eq:reduced_set}
\end{equation}
The feature block associated with $j$ is the submatrix containing the linear term for $j$, its square, and every pairwise interaction involving $j$,
\begin{equation}
\mathbf{G}_m\!\left(j;\mathcal{B}^{(q-1)}_m\right)
=
\left[
\overline{X}_{m,j}(t),
\overline{X}_{m,j}^{2}(t),
\overline{X}_{m,j}(t)\overline{X}_{m,k}(t)
\;\middle|\;
 k\in\mathcal{B}^{(q-1)}_{m,-j}
\right]_{t=1}^{T-1}.
\label{eq:variable_block}
\end{equation}
This block has $|\mathcal{B}^{(q-1)}_m|+1$ columns. The conditional contribution assigned to $j$ is
\begin{equation}
\Delta I^{-,(q)}_{m,j}
=
I\!\left(
\mathbf{y}_m;
\mathcal{P}_m\!\left(\mathbf{G}_m(j;\mathcal{B}^{(q-1)}_m)\right)
\,\middle|\,
\mathcal{P}_m\!\left(\mathbf{Q}_m(\mathcal{B}^{(q-1)}_{m,-j})\right)
\right).
\label{eq:deltaI_backward}
\end{equation}
The least contributing variable is removed,
\begin{equation}
j_q^{-}
=
\arg\min_{j\in\mathcal{B}^{(q-1)}_m}
\Delta I^{-,(q)}_{m,j},
\qquad
\mathcal{B}^{(q)}_m
=
\mathcal{B}^{(q-1)}_m\setminus\{j_q^{-}\}.
\label{eq:update_backward}
\end{equation}
The quadratic design is rebuilt after each removal. Two sequential backward iterations are used, so $K_b=2$, and the final local variable set is
\begin{equation}
\mathcal{S}_m
=
\mathcal{B}^{(K_b)}_m,
\qquad
|\mathcal{S}_m|=K_f-K_b=8.
\label{eq:final_local_set}
\end{equation}
The final eight-variable expansion contains 44 linear and quadratic terms. The backward stage therefore validates variable groups under the nonlinear representation while preserving a physically interpretable output at the variable level.

\subsection{Cardinality choice}
\label{sec:cardinality_constraint}

The original ER formulation uses an information-based tolerance to stop the forward and backward stages \cite{almomani2020entropic,almomani2020erfit}. The present comparison instead fixes the final cardinality so that every superpixel and every method returns the same number of variables. This design makes the recurrence profiles directly comparable, but it does not imply that every local region has the same intrinsic structural dimension.

For each superpixel, form the predictor matrix
\begin{equation}
\mathbf{X}_m\in\mathbb{R}^{T\times p},
\qquad
(\mathbf{X}_m)_{t,j}=\overline{X}_{m,j}(t),
\label{eq:Xm}
\end{equation}
and compute
\begin{equation}
\mathbf{X}_m
=
\mathbf{U}_m\boldsymbol{\Sigma}_m\mathbf{V}_m^{\top}.
\label{eq:svd_X}
\end{equation}
The cumulative singular-value energy is
\begin{equation}
E_m(k)
=
\frac{\sum_{i=1}^{k}\sigma_{m,i}^2}
{\sum_{i=1}^{r_m}\sigma_{m,i}^2},
\label{eq:energy_local}
\end{equation}
and its domain average is
\begin{equation}
E(k)
=
\frac{1}{M}\sum_{m=1}^{M}E_m(k).
\label{eq:energy_global}
\end{equation}
The retained analysis gave
\begin{equation}
E(8)\approx0.992.
\label{eq:energy_result}
\end{equation}
Accordingly, $K=8$ was adopted as a working low-rank cardinality. Principal component and singular-value analyses describe dominant variance directions, not causal cardinality \cite{Jolliffe2016PCA}. The calculation is therefore used only as a complexity heuristic. It neither establishes that eight original variables are uniquely sufficient nor excludes low-variance fields from structural relevance.

\subsection{Information estimation and scope of the score}
\label{sec:information_estimation}

The lagged samples are constructed from $\{(\overline{\mathbf{X}}_m(t),\overline{Y}_m(t+1))\}_{t=1}^{T-1}$. Each predictor series is standardized within its superpixel using the development sample. Mutual information and conditional mutual information are estimated with a nearest-neighbor estimator following the Kraskov, St\"ogbauer, and Grassberger framework and its conditional extension \cite{kraskov2004estimating,frenzel2007partial}. The estimator is applied to scalar regression projections and the scalar target, rather than directly to the full predictor set.

This reduction in estimator dimension is the principal computational motivation for ER. It is also important to state what the score is not. In general,
\begin{equation}
I\!\left(
\mathbf{y}_m;
\mathbf{z}_{m,j}
\,\middle|\,
\mathbf{z}_{m,\mathcal{S}}
\right)
\neq
I\!\left(
\overline{X}_{m,j};
\overline{Y}_m(t+1)
\,\middle|\,
\overline{\mathbf{X}}_{m,\mathcal{S}}
\right).
\label{eq:projection_not_raw}
\end{equation}
The left side is a projection-mediated causal-relevance score. It asks whether a candidate contributes unique information through the adopted model class after conditioning on the current selected projection. The resulting causal structure is therefore defined relative to the feature representation and projection rule, rather than by raw-variable conditional independence alone.

Nearest-neighbor information estimation also requires caution in the present application. Precipitation contains many exact zeros and therefore combines an atom at zero with a continuous positive component. Estimation for mixed and repeated-value data is a distinct problem from estimation under a smooth continuous density \cite{Gao2017Mixture,Mesner2021MixedCMI}. The geometric entropy framework and Boltzmann-Shannon interaction entropy provide related alternatives for sparse samples, informative extremes, and partition-based information analysis \cite{diggans2022geometric,diggans2023boltzmann,diggans2025generalizing}, but those estimators were not substituted into the retained experiment. The reported scores should consequently be understood as outputs of the implemented nearest-neighbor procedure.

Finally, the regression projections and information scores were estimated on the same development sample. No blocked cross-fitting was used. This preserves the retained implementation but permits in-sample optimism in the information scores. Temporally blocked cross-fitting is an important requirement for a future confirmatory study because random splitting is generally inadequate for dependent spatial and temporal data \cite{Roberts2017CV}.

\subsection{Domain-level recurrence}
\label{sec:domain_level_aggregation}

Applying the procedure independently to the $M$ superpixels produces the sets $\{\mathcal{S}_m\}_{m=1}^{M}$. For variable $j$, define the selection proportion
\begin{equation}
\pi_j
=
\frac{1}{M}
\sum_{m=1}^{M}
\mathbf{1}_{\{j\in\mathcal{S}_m\}}.
\label{eq:selection_frequency}
\end{equation}
The statistic $\pi_j$ is the fraction of superpixels in which variable $j$ is selected. A large value identifies a candidate parent whose conditional contribution is not confined to isolated regions but persists across the domain. Its interpretation is related in spirit to stability-selection ideas that summarize recurrence across perturbed data sets \cite{Meinshausen2010Stability}. The present construction is not formal stability selection, however. The superpixels are spatially dependent, they arise from one fixed precipitation-guided partition, and no subsampling-based false-selection control is claimed. Thus, $\pi_j$ measures spatial recurrence of the recovered causal structure rather than a calibrated probability of causal truth.

Variables with large $\pi_j$ define the domain-level causal backbone recovered by ER. Variables with smaller $\pi_j$ may be localized, redundant with recurrent fields, sensitive to the local data distribution, or weak under the chosen model class. The recurrence statistic alone does not distinguish among these possibilities.

\section{Comparison with Transfer Entropy and Causation Entropy}
\label{sec:comparison}

The comparison tests whether three information-theoretic procedures recover a concentrated and spatially recurrent candidate-parent structure under the same local cardinality. Because no interventionally verified precipitation graph is available, the comparison evaluates the geometry and recurrence of the recovered causal structures rather than absolute graph-recovery accuracy.

\subsection{Transfer entropy}
\label{sec:transfer_entropy}

For stochastic processes $X$ and $Y$, transfer entropy is
\begin{equation}
T_{X\rightarrow Y}
=
I\!\left(X_t^{-};Y_{t+1}\mid Y_t^{-}\right),
\label{eq:TE_def}
\end{equation}
where $X_t^{-}$ and $Y_t^{-}$ denote process histories. Under a first-order representation,
\begin{equation}
T_{X\rightarrow Y}
=
H(Y_{t+1}\mid Y_t)
-
H(Y_{t+1}\mid Y_t,X_t).
\label{eq:TE_markov}
\end{equation}
Transfer entropy measures directed predictive information beyond target history \cite{Schreiber2000}. For Gaussian processes it is closely related to Granger causality \cite{Barnett2009}, but no Gaussian assumption is imposed here. In its pairwise form, a predictor can receive a large score because of indirect pathways, shared drivers, or redundant representations of the same atmospheric state.

\subsection{Causation entropy}
\label{sec:causation_entropy}

For index sets $I$, $J$, and $K$, causation entropy is
\begin{equation}
C_{J\rightarrow I\mid K}
=
H\!\left(X_{t+1}^{(I)}\mid X_t^{(K)}\right)
-
H\!\left(X_{t+1}^{(I)}\mid X_t^{(K)},X_t^{(J)}\right),
\label{eq:CE_def}
\end{equation}
or equivalently,
\begin{equation}
C_{J\rightarrow I\mid K}
=
I\!\left(X_t^{(J)};X_{t+1}^{(I)}\mid X_t^{(K)}\right).
\label{eq:CE_cmi}
\end{equation}
The optimal causation entropy principle identifies a minimal conditioning set that maximizes the relevant information quantity under its theoretical assumptions \cite{Sun2015OCSE}. The same forward and backward logic has been adapted to data-driven Boolean-network and Boolean-function inference \cite{sun2022data}. In continuous high-dimensional applications, the conditioning dimension grows as variables are admitted, which increases the finite-sample burden on the conditional-information estimator \cite{frenzel2007partial,Runge2018CMI}.

\subsection{Relation to entropic regression}
\label{sec:relation_to_er}

ER and CE share the principle of iterative conditional-information screening, but their operational quantities differ. CE estimates conditional information from the candidate variable and the current conditioning set. ER first maps a candidate term or variable group to a fitted response projection and then evaluates the information carried by that projection. Consequently, ER is not an algebraically equivalent low-dimensional form of raw-variable causation entropy. It is a regression-mediated screening procedure motivated by the same incremental-information principle \cite{almomani2020entropic,fish2021entropic}.

\subsection{Comparison protocol}
\label{sec:comparison_protocol}

Each method is applied to the same collection of superpixel-level time series and the same one-hour predictor-response pairing. Every local result is restricted to
\begin{equation}
|\mathcal{S}_m|=K=8.
\label{eq:protocol_cardinality}
\end{equation}
Transfer entropy retains its target-history conditioning, while CE and ER use their respective implemented conditioning sets. Thus, the common cardinality controls selected-set size but does not make the estimands identical. For each method, Eq.~\eqref{eq:selection_frequency} is used to compute a domain-level selection proportion. The comparison is therefore descriptive: it asks how concentrated and recurrent the outputs are under the implemented procedures, not which method recovers an unknown true graph.

\section{Results}
\label{sec:results}

\subsection{Spatial recurrence of candidate causal parents across methods}
\label{sec:results_structural_consistency}

\begin{figure}[t]
\centering
\includegraphics[width=0.94\columnwidth]{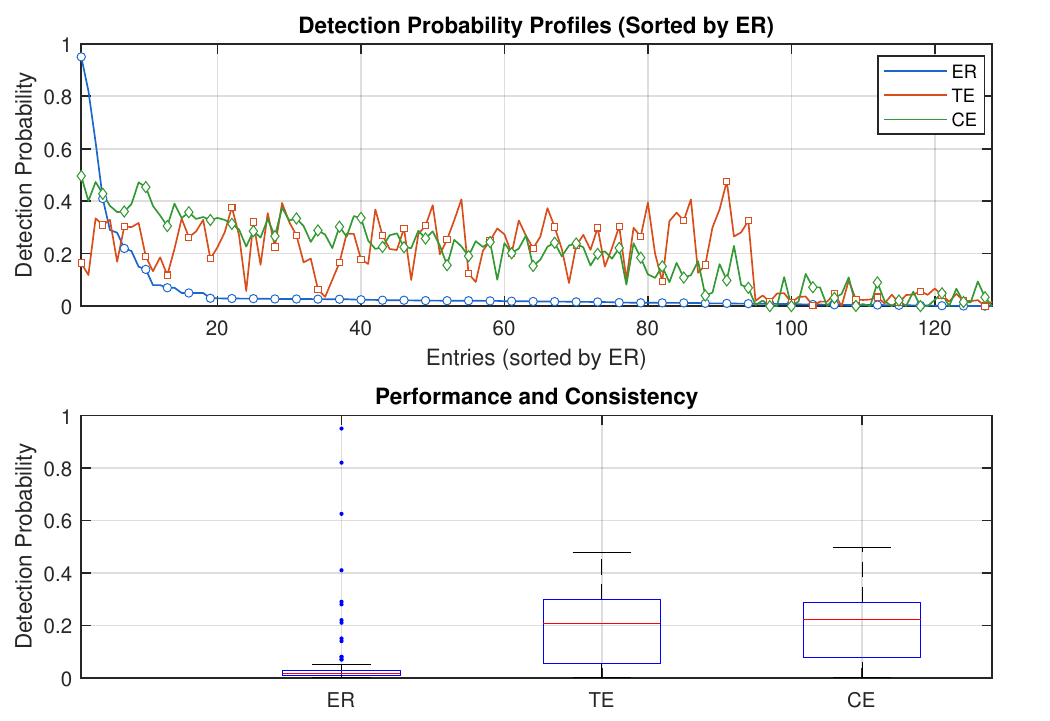}
\caption{Domain-level selection proportions for entropic regression, transfer entropy, and causation entropy. Top: variable-wise profiles ordered by decreasing ER recurrence. Bottom: boxplot summaries of the three recurrence distributions. All methods return eight variables per superpixel.}
\label{fig:method_comparison}
\end{figure}

Figure~\ref{fig:method_comparison} compares the selection proportions obtained from ER, TE, and CE. Variables are ordered by decreasing ER recurrence, and the same ordering is used for the other methods.

The ER profile is sharply concentrated. A small number of variables recur in a large fraction of the superpixels, followed by a rapid decrease toward a broad group with low recurrence. Because every method returns eight variables in every superpixel, this concentration is not a consequence of a smaller local model. It shows that ER repeatedly resolves the same domain-level parent structure across heterogeneous local inference problems.

The TE profile is substantially more diffuse across the ordered variables. CE also distributes its selections across a wider set than ER, although its ranking and dispersion differ from TE. Under the implemented protocol, ER therefore produces the clearest separation between recurrent candidate parents and context-dependent variables.

This comparative result is positive and specific: projection-based ER converts a high-dimensional field set into a compact causal-parent recurrence profile, while the direct-information alternatives retain broader domain-level support. The absence of a known true precipitation graph prevents an absolute accuracy ranking, but it does not alter the observed difference in structural concentration.

The principal comparative result is that, under the same local cardinality, ER recovers a substantially more concentrated and recurrent parent structure. This provides a direct answer to the structural question posed in the paper: projection-based ER resolves a compact domain-level causal backbone from a 122-variable atmospheric state. The backbone is summarized in the following section.

\subsection{Recurrent structural variables}
\label{sec:results_dominant_drivers}

To extract a compact domain-level summary, define the operational recurrence threshold
\begin{equation}
\pi_j\geq\tau,
\qquad
\tau=0.25.
\label{eq:threshold}
\end{equation}
Thus, a reported variable must be selected in at least one quarter of the superpixels. The threshold is descriptive. It is not a universal significance level or a causal probability cutoff.

Six variables satisfy Eq.~\eqref{eq:threshold}. These six fields define the recurrent causal structure recovered by ER at the domain level. Each has survived sequential conditional-information selection within local models and recurs in at least one quarter of the 1227 superpixels, making the result stronger than a marginal correlation ranking. Their retained field labels and recurrence proportions are shown in Table~\ref{tab:causal_parents}. The descriptions follow the retained HRRR field labels and the official product inventory \cite{dowell2022hrrr,NOAAHRRRInventory}.

\begin{table*}[t]
\centering
\caption{Variables exceeding the operational recurrence threshold $\pi_j\geq0.25$. The selection proportion is the fraction of the 1227 superpixels in which the variable appears in the final eight-variable ER set.}
\label{tab:causal_parents}
\begin{tabular}{l p{9.2cm} l c}
\toprule
\textbf{Field} & \textbf{Description} & \textbf{Unit} & \textbf{$\pi_j$} \\
\midrule
MSTAV & Land-surface moisture availability & \% & 0.95 \\
REFC & Composite reflectivity & dB & 0.82 \\
DZDT & Geometric vertical velocity & m/s & 0.62 \\
WIND & Wind speed & m/s & 0.41 \\
CAPE & Convective available potential energy & J/kg & 0.29 \\
MAXUVV & Hourly maximum upward vertical velocity in the lowest 400 hPa & m/s & 0.28 \\
\bottomrule
\end{tabular}
\end{table*}

The six fields are physically interpretable, but their roles are not equivalent. REFC is composite reflectivity, an active-state indicator of hydrometeor and precipitation organization. Because HRRR assimilates radar information, REFC may carry strong persistence and analysis-state information rather than acting as an upstream environmental driver in the same sense as instability or vertical motion \cite{dowell2022hrrr,NOAAHRRRInventory}. MSTAV is the modeled land-surface moisture-availability field, not atmospheric water-vapor content. Its recurrence may reflect land-atmosphere coupling, antecedent wetness, spatial climatology, or memory associated with prior precipitation \cite{Koster2004,Seneviratne2010,Ford2015}. DZDT and MAXUVV represent aspects of vertical motion, WIND is the retained wind-speed field, and CAPE characterizes parcel instability. These quantities are physically consistent with lifting, flow organization, and convective development, but their precise role depends on level, regime, and the surrounding atmospheric state \cite{Houze2014,Markowski2010}.

Taken together, the variables span land-surface state, active hydrometeor structure, vertical motion, wind speed, and thermodynamic instability. Their recurrence is physically coherent with precipitation-relevant conditions, but the analysis does not isolate mechanisms or establish that any field is individually necessary. Correlated substitutes and unresolved common drivers remain possible.

The six-variable set is therefore the \emph{recurrent causal backbone} recovered by ER for the adopted precipitation representation. It is the compact parent structure that repeatedly survives sequential conditional-information selection across the domain. The term does not require that the six variables exhaust every physical pathway or unresolved state variable.

\subsection{Predictive consequences of the recovered causal structure}
\label{sec:predictive_outcome}

The recurrent causal backbone is next used as a reduced predictor representation. This analysis tests the predictive consequences of the discovered structure: which aspects of next-hour precipitation remain accessible when the candidate causal parents are used as the only inputs to the tested models.

\subsubsection{Precipitation occurrence}
\label{sec:probability_of_precipitation}

Define the binary target
\begin{equation}
Z(\mathbf{s},t)
=
\mathbf{1}_{\{Y(\mathbf{s},t)>0\}}.
\label{eq:binary_target}
\end{equation}
For forecast horizon $\tau$, the model estimates
\begin{equation}
\widehat{p}_{\tau}(\mathbf{s})
=
\mathbb{P}\!\left(
Z(\mathbf{s},t+\tau)=1
\mid
\text{selected predictor history}
\right).
\label{eq:pop_model}
\end{equation}
A convolutional sequence model is used to combine localized spatial features with temporal history. The recurrent component follows the LSTM family of sequence models \cite{hochreiter1997lstm}, and convolutional recurrence has been widely used for precipitation nowcasting \cite{shi2015convlstm,ayzel2020rainnet,ravuri2021skilful}.

\begin{figure}[t]
\centering
\includegraphics[width=0.88\columnwidth]{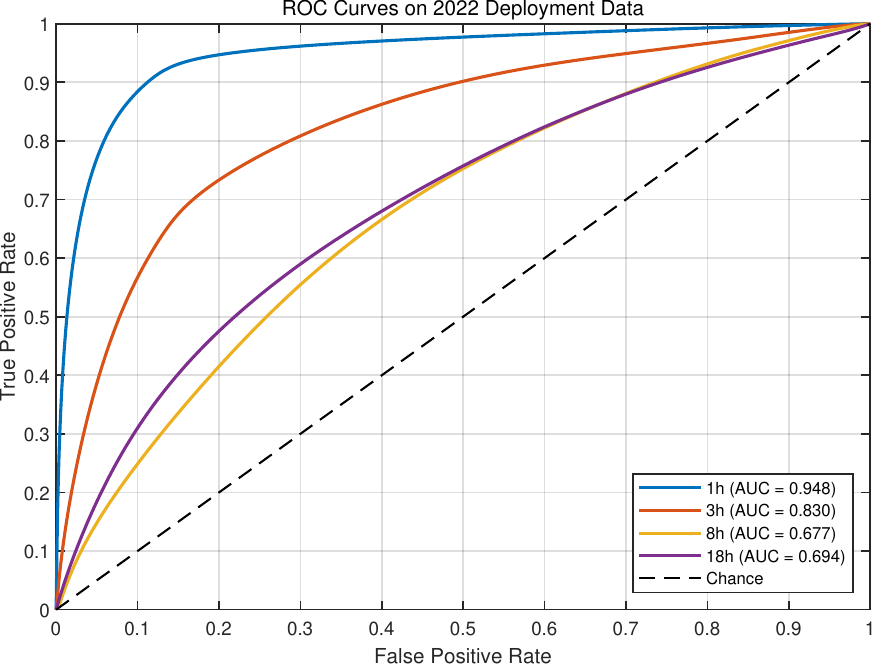}
\caption{Receiver operating characteristic curves for precipitation occurrence on the 2022 deployment sample. The area under the curve summarizes threshold-independent discrimination for forecast horizons of 1, 3, 8, and 18 hours.}
\label{fig:roc_probability_model}
\end{figure}

Figure~\ref{fig:roc_probability_model} reports receiver operating characteristic curves for $\tau\in\{1,3,8,18\}$ hours on the 2022 deployment sample. Receiver operating characteristic analysis summarizes ranking discrimination over decision thresholds \cite{Fawcett2006}. The one-hour area under the curve is 0.948, and the three-hour value is 0.830. The corresponding values at eight and eighteen hours are 0.677 and 0.694. Thus, the recovered causal structure retains substantial short-horizon event information, with particularly strong one-hour discrimination, while its ability to rank future precipitation states decreases at the longer horizons considered.

The result should be interpreted as event discrimination rather than full probabilistic sufficiency. A high area under the curve indicates that precipitation and non-precipitation cases are ranked effectively. It does not establish calibrated probability values, prevalence-sensitive precision-recall behavior, or accurate precipitation amounts \cite{Murphy1993,Gneiting2007,Saito2015}.

\begin{figure}[p]
\centering
\includegraphics[width=0.95\textwidth,height=0.70\textheight,keepaspectratio]{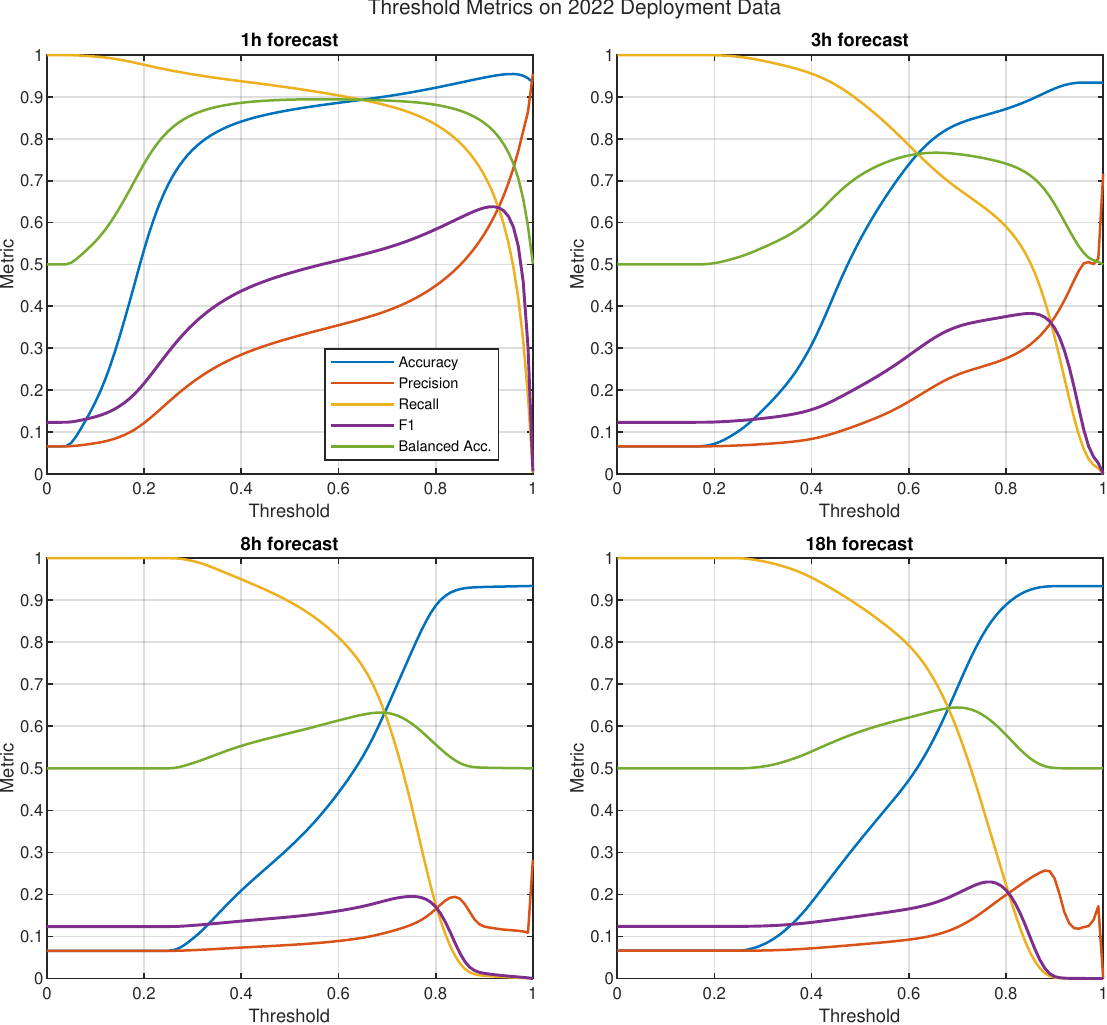}
\caption{Classification metrics as functions of the probability threshold for the four forecast horizons. The variation across thresholds emphasizes that discrimination, precision, recall, and balanced accuracy describe different properties of the same probabilistic output.}
\label{fig:threshold_metrics}
\end{figure}

Figure~\ref{fig:threshold_metrics} shows the associated threshold dependence of accuracy, precision, recall, F1 score, and balanced accuracy. Such metrics emphasize different error tradeoffs and are especially sensitive to event prevalence and operating threshold \cite{Saito2015,Wilks2019}. The curves demonstrate that the apparent performance depends materially on the operating threshold, especially under the strong event imbalance. This threshold sensitivity is part of the predictive characterization rather than evidence for a single preferred deployment rule.

\subsubsection{Intensity and spatial summaries}
\label{sec:intensity_spatial}

A convolutional LSTM model was also used to predict the precipitation field at a one-hour horizon from the six selected variables. Let
\begin{equation}
\widehat{Y}(t+1)
=
\mathcal{G}\!\left(\{\mathbf{X}_{\mathcal{S}}(t-\ell)\}_{\ell=0}^{T_h}\right),
\label{eq:lstm_model}
\end{equation}
where $\mathcal{G}$ denotes the fitted spatiotemporal model and $\mathbf{X}_{\mathcal{S}}$ contains only the selected fields.

Over the full deployment sample, the retained one-hour intensity metrics are
\begin{equation}
\mathrm{RMSE}=0.8528,
\quad
\mathrm{MAE}=0.1014,
\quad
\mathrm{Bias}=-0.0308,
\quad
\rho=0.2894,
\quad
R^2=0.0574.
\label{eq:intensity_metrics}
\end{equation}
When evaluation is restricted to observations with rain, the retained metrics are
\begin{equation}
\mathrm{RMSE}=3.3098,
\quad
\mathrm{MAE}=1.0743,
\quad
\rho=0.2038.
\label{eq:rain_only_metrics}
\end{equation}
The larger rain-only errors show that much of the full-sample performance is associated with the dominant zero and low-precipitation regime. The selected representation retains some information about where precipitation occurs, but the tested model resolves only a limited fraction of conditional intensity variation.

The domain-level spatial summaries are
\begin{equation}
\rho_{\mathrm{spatial}}=0.7781,
\quad
\mathrm{RMSE}_{\mathrm{spatial}}=0.0558,
\quad
\mathrm{MAE}_{\mathrm{spatial}}=0.0373,
\quad
\mathrm{Bias}_{\mathrm{spatial}}=-0.0301.
\label{eq:spatial_metrics}
\end{equation}
The relatively high spatial correlation, together with weak pointwise intensity correlation, indicates that broad spatial organization is represented more successfully than local amplitude. Spatial precipitation verification is scale dependent, and a field may reproduce broad organization while missing local displacement or amplitude \cite{Gilleland2009}. The exact meaning of this contrast depends on the aggregation used to compute Eq.~\eqref{eq:spatial_metrics}, so these metrics are interpreted as complementary summaries rather than as interchangeable measures of forecast quality.

Overall, the predictive experiments reveal a differentiated but scientifically informative profile. The causal backbone retains strong one-hour occurrence discrimination and broad spatial information, while intensity resolution under the tested architecture is much weaker. This establishes the second central result of the paper: causal relevance and predictive sufficiency are distinct properties. The discovered parents identify a common information-bearing structure of the precipitation process, while accurate intensity evolution requires additional conditional structure, temporal memory, unresolved variables, or a different predictive representation.

\section{Discussion}
\label{sec:discussion}

\subsection{What the causal discovery result establishes}

The paper establishes more than a marginal association pattern. Within every superpixel, ER admits a variable only when its projection contributes information beyond the projection formed by the variables already selected. Six variables then recur in at least one quarter of 1227 spatially heterogeneous local problems. Sequential conditional contribution and spatial recurrence therefore identify a compact, domain-level causal backbone. This is the central causal discovery result of the study.

This result provides a scientifically useful reduction of the 122-variable atmospheric state. It identifies the fields that organize the one-hour precipitation response most persistently under the adopted model class, rather than merely listing predictors with large pairwise dependence. The recurrence profile also gives a concrete set of physical hypotheses for subsequent analysis of precipitation initiation, persistence, and intensity.

The distinction between recurrence and calibrated causal probability remains important. Formal stability selection obtains error-control properties by applying a selection procedure to repeated random subsamples under explicit assumptions \cite{Meinshausen2010Stability}. The present superpixels are neither random subsamples nor independent perturbations. They are spatially dependent regions from one data-guided partition. Consequently, $\pi_j$ is a spatial recurrence proportion, not a posterior probability that a variable is causally true. This boundary limits the statistical interpretation of the number, but it does not reduce the substantive finding that the same parent variables repeatedly emerge across the domain.

\subsection{Physical interpretation and alternative pathways}

The recurrent causal backbone is physically coherent, while containing fields that occupy different positions within the precipitation process. Vertical velocity, maximum upward vertical velocity, wind speed, and CAPE are associated with lifting, circulation, and instability that organize convection \cite{Houze2014,Markowski2010}. This coherence supports physical plausibility, but it does not prove that the ER score has isolated a mechanism. Correlated fields can substitute for one another, and a selected field can be an effect, proxy, or state indicator rather than an upstream cause.

REFC requires particular care. Composite reflectivity describes current hydrometeor organization, and HRRR analyses incorporate radar observations \cite{dowell2022hrrr,NOAAHRRRInventory}. Its high recurrence is therefore consistent with strong one-hour persistence information. That result is useful for structural screening, but REFC should not be interpreted in the same causal category as an environmental precursor without additional conditioning on current precipitation or reflectivity.

MSTAV also admits several pathways. The field represents land-surface moisture availability. Soil moisture can influence surface fluxes, boundary-layer development, and convective initiation, and the strength and sign of soil-moisture precipitation coupling vary across region and regime \cite{Koster2004,Seneviratne2010,Ford2015}. At the same time, surface wetness can retain memory of earlier precipitation. The present one-hour observational design cannot distinguish land-atmosphere feedback from antecedent-rainfall memory or spatial climatology.

This ambiguity does not make the recurrence result uninformative. It suggests a layered hypothesis: the recurrent set may combine active-state indicators, environmental conditions, and memory variables. A future study can test that hypothesis by conditioning on target history and current reflectivity, excluding diagnostic precipitation fields, and then screening the residual information carried by the remaining atmospheric variables.

\subsection{Causal relevance and predictive sufficiency}

The second principal contribution is the empirical separation between causal relevance and predictive sufficiency. The recovered parents retain strong one-hour occurrence discrimination, yet the tested model resolves only a limited fraction of pointwise intensity variation. These results show that a variable can be causally informative under a lagged conditional-information criterion without forming a complete predictive state for every forecast quantity.

This is not a negative result. It identifies two layers of the precipitation problem. The recurrent causal backbone captures a common domain-level structure that distinguishes active and inactive precipitation states and organizes broad spatial patterns. Additional conditional structure becomes important within active events, where precipitation amount depends on storm evolution, microphysics, terrain, unresolved scales, temporal history, and regime-specific interactions.

The broader forecast-verification literature likewise distinguishes discrimination, reliability, threshold utility, quantitative accuracy, and spatial structure as different attributes \cite{Murphy1993,Gneiting2007,Fawcett2006,Saito2015,Gilleland2009,Wilks2019}. The present paper connects that distinction to causal discovery. A reduced causal representation should be evaluated by asking which predictive attributes it retains, rather than by demanding that one structural criterion simultaneously optimize every forecast metric.

The current experiments do not isolate which missing component is responsible for the weaker intensity result. Plausible explanations include missing variables, insufficient temporal history, imperfect optimization, unresolved microphysics, spatial misalignment, target uncertainty, class imbalance, and limits of the chosen architecture. Without matched all-variable, persistence-only, reflectivity-only, and random-subset baselines, these explanations cannot be separated.

A constructive hypothesis follows directly from the results: the six recurrent parents form a common causal layer, while the predictive residual contains event-conditioned structure. Applying causal discovery to rain-only intensity residuals, or conditioning explicitly on current precipitation and reflectivity, may reveal distinct variables associated with initiation, maintenance, and intensity. The present paper therefore provides both a causal backbone and a precise target for the next study.

\subsection{Comparison with transfer entropy and causation entropy}

Under the common-cardinality protocol, the concentrated ER profile is a practical methodological advantage. ER resolves a small set of parents that recurs across the domain, whereas TE and CE distribute their local selections across broader variable sets. Pairwise TE can assign information to several related fields that share indirect or common-driver pathways \cite{Schreiber2000}. CE introduces multivariate conditioning and is theoretically designed to separate direct from indirect influence under its assumptions \cite{Sun2015OCSE}, but its empirical burden increases with the dimension of the conditioning space \cite{frenzel2007partial,Runge2018CMI}. ER instead compresses candidate information through fitted projections \cite{almomani2020entropic,fish2021entropic}.

The concentrated ER profile may therefore arise because regression projection compresses related variables into a lower-dimensional fitted representation. This is the intended computational advantage, but it is also a modeling commitment. If the linear forward model or quadratic backward library fails to represent an important dependence, the corresponding variable may receive a small ER score even when it is informative under another representation. Concentration should not be equated automatically with correctness.

The supported comparative claim is therefore clear: in this high-dimensional precipitation problem, projection-based ER recovers the most concentrated and spatially recurrent candidate-parent structure among the three implemented methods. Establishing universal superiority or absolute graph-recovery accuracy would require controlled synthetic benchmarks, alternate partitions, repeated samples, and known ground truth.

\subsection{Limitations and future directions}
\label{sec:limitations}

Several limitations define the scope of the results.

First, the structural interpretation is observational and model dependent. The ER score is evaluated through fitted projections and is not identical to raw-variable conditional mutual information. No intervention or ground-truth precipitation graph is available \cite{Pearl2009,Peters2017,Runge2019CausationEarth}.

Second, the forward projections and their information scores are estimated from the same development sample. This can introduce in-sample optimism. A confirmatory study should use temporally blocked cross-fitting, with projection coefficients estimated outside the blocks used for information evaluation \cite{Roberts2017CV}.

Third, the final cardinality is fixed at eight. The SVD analysis provides a low-rank heuristic, but variance concentration does not determine causal cardinality. Low-variance variables may be relevant, and different regions may require different model sizes \cite{Jolliffe2016PCA}.

Fourth, the forward stage uses only linear variable projections, and the nonlinear model is restricted to quadratic terms after forward screening. A variable whose relevance appears only through a nonlinear interaction may therefore be missed before the expansion is constructed. Conversely, the adaptive design avoids a 7625-term global library and preserves variable-level interpretability. This tradeoff should be tested directly in future work.

Fifth, the recurrence proportions are obtained from one precipitation-guided partition. The superpixels are spatially dependent, and no alternate partition, spatial block bootstrap, or formal stability-selection analysis is available. The results establish recurrence across the adopted regions, not stability to arbitrary spatial representations.

Sixth, nearest-neighbor information estimates are sensitive to finite samples, temporal dependence, repeated values, and estimator settings. Precipitation contains a point mass at zero, so the numerical scores are estimator dependent. Mixed-data estimators and geometry-aware alternatives provide useful directions for sensitivity analysis \cite{Gao2017Mixture,Mesner2021MixedCMI,diggans2022geometric,diggans2023boltzmann,diggans2025generalizing}.

Seventh, the predictive experiments do not include matched models using all 122 variables, persistence-only inputs, reflectivity alone, or randomly selected subsets. The results characterize the six-variable representation under the tested models but do not identify the source of remaining error.

These limitations motivate a direct follow-up program. Future work should examine alternate partitions, cardinality sensitivity, explicit target-history conditioning, reflectivity-excluded screening, temporally blocked cross-fitting, lag-resolved selection, and event-conditioned analysis. Applying structural screening to rain-only intensity residuals may reveal variables hidden by the dominant occurrence process. Such analyses can test whether the present domain-level profile is supplemented by distinct initiation, maintenance, and intensity structures.

\section{Conclusion}
\label{sec:conclusion}

This study makes two primary contributions. First, it introduces a projection-based, variable-level formulation of entropic regression for structural causal discovery in high-dimensional nonlinear systems. The method performs ten sequential forward selections on the physical variables, constructs a quadratic library only over the resulting reduced set, and applies two sequential variable-block eliminations. It therefore preserves the forward and backward causal-selection logic of entropic regression while avoiding construction of the 7625-term global quadratic library.

Second, the method recovers a compact and spatially recurrent causal structure in precipitation dynamics. Across 1227 precipitation-guided superpixels, six fields exceed the operational recurrence threshold and form a physically coherent backbone involving land-surface state, hydrometeor organization, vertical motion, wind speed, and convective instability. Under the same final cardinality, transfer entropy and causation entropy produce broader recurrence profiles. The result shows that projection-based ER can resolve an interpretable domain-level parent structure from a strongly coupled 122-variable atmospheric state.

The predictive analysis establishes the conceptual contribution of the paper. The recovered causal backbone retains strong one-hour precipitation-occurrence discrimination and broad spatial information, while pointwise intensity prediction remains substantially more limited under the tested architecture. Causal relevance and predictive sufficiency are therefore distinct properties of a reduced dynamical representation. The first identifies variables that contribute recurrent, nonredundant information about the next state. The second requires enough state information to reproduce the particular forecast quantity and scale of interest.

This distinction gives the structural result practical scientific value. The six-variable backbone is not merely a compressed forecast input set. It is a testable hypothesis about the common information layer organizing next-hour precipitation across the Southwestern United States. The unresolved predictive variation points to additional event-conditioned structure, including target history, diagnostic persistence, regime dependence, terrain, microphysics, and unresolved scales. The present study therefore closes one stage of the causal analysis while defining the next: separating the common causal backbone from the conditional structures governing initiation, maintenance, and intensity.

The term candidate causal parent retains the operational observational meaning defined in this paper. Within that scope, the contribution is affirmative and specific: variable-level entropic regression identifies a compact recurrent causal backbone in a high-dimensional Earth-system data set. The predictive analysis then shows that this backbone captures a common causal layer for precipitation occurrence, while additional conditional structure is required to resolve intensity evolution. Structural discovery and predictive modeling are therefore complementary, rather than interchangeable, scientific tasks.

\subsection*{Data and Code Availability}
\label{sec:data_code_availability}

The HRRR atmospheric fields and MRMS precipitation observations are publicly available from their respective data providers. The project repository contains the current manuscript, the retained predictor-variable list, and workflow documentation:
\begin{center}
\url{https://github.com/almomaa/EAGER_ER}
\end{center}
The general ERFit MATLAB implementation is maintained separately at:
\begin{center}
\url{https://github.com/almomaa/ERFit-Package}
\end{center}
Large processed atmospheric and precipitation fields are not included in the manuscript repository. The complete original execution environment is not claimed to be archived. Retained outputs and additional project materials may be provided by the lead author subject to availability.

\section*{Acknowledgments}

This material is based upon work supported by the National Science Foundation under Grant No. 2313689. Any opinions, findings, conclusions, or recommendations expressed in this material are those of the authors and do not necessarily reflect the views of the National Science Foundation.

\section*{Declaration of Interests}

The authors declare no competing interests.


\bibliographystyle{unsrt}
\bibliography{references}

@book{Kalnay2003,
  title     = {Atmospheric Modeling, Data Assimilation and Predictability},
  author    = {Kalnay, Eugenia},
  year      = {2003},
  publisher = {Cambridge University Press},
  doi       = {10.1017/CBO9780511802270}
}

@article{dowell2022hrrr,
  title   = {The High-Resolution Rapid Refresh (HRRR): An hourly updating convection-allowing forecast model. Part I: Motivation and system description},
  author  = {Dowell, David C. and Alexander, Curtis R. and James, Eric P. and Weygandt, Stephen S. and Benjamin, Stan G. and Manikin, Geoff S. and Blake, Benjamin T. and Brown, John M. and Olson, Joseph B. and Hu, Ming and others},
  journal = {Weather and Forecasting},
  volume  = {37},
  number  = {8},
  pages   = {1371-1395},
  year    = {2022},
  doi     = {10.1175/WAF-D-21-0151.1}
}

@article{reichstein2019deep,
  title   = {Deep learning and process understanding for data-driven Earth system science},
  author  = {Reichstein, Markus and Camps-Valls, Gustau and Stevens, Bjorn and Jung, Martin and Denzler, Joachim and Carvalhais, Nuno and Prabhat},
  journal = {Nature},
  volume  = {566},
  pages   = {195-204},
  year    = {2019},
  doi     = {10.1038/s41586-019-0912-1}
}

@article{rasp2020weatherbench,
  title   = {WeatherBench: A benchmark data set for data-driven weather forecasting},
  author  = {Rasp, Stephan and Dueben, Peter D. and Scher, Sebastian and Weyn, Jonathan A. and Mouatadid, Soukayna and Thuerey, Nils},
  journal = {Journal of Advances in Modeling Earth Systems},
  volume  = {12},
  number  = {11},
  pages   = {e2020MS002203},
  year    = {2020},
  doi     = {10.1029/2020MS002203}
}

@article{rasp2024weatherbench2,
  title   = {WeatherBench 2: A benchmark for the next generation of data-driven global weather models},
  author  = {Rasp, Stephan and others},
  journal = {Journal of Advances in Modeling Earth Systems},
  volume  = {16},
  number  = {6},
  pages   = {e2023MS004019},
  year    = {2024},
  doi     = {10.1029/2023MS004019}
}

@article{ayzel2020rainnet,
  title   = {RainNet v1.0: A convolutional neural network for radar-based precipitation nowcasting},
  author  = {Ayzel, Georgy and Scheffer, Tobias and Heistermann, Maik},
  journal = {Geoscientific Model Development},
  volume  = {13},
  number  = {6},
  pages   = {2631-2644},
  year    = {2020},
  doi     = {10.5194/gmd-13-2631-2020}
}

@inproceedings{shi2015convlstm,
  title     = {Convolutional LSTM Network: A Machine Learning Approach for Precipitation Nowcasting},
  author    = {Shi, Xingjian and Chen, Zhourong and Wang, Hao and Yeung, Dit-Yan and Wong, Wai-kin and Woo, Wang-chun},
  booktitle = {Advances in Neural Information Processing Systems},
  volume    = {28},
  pages     = {802-810},
  year      = {2015}
}

@article{ravuri2021skilful,
  title   = {Skilful precipitation nowcasting using deep generative models of radar},
  author  = {Ravuri, Suman and Lenc, Karel and Willson, Matthew and Kangin, Dmitry and Lam, Remi and Mirowski, Piotr and Fitzsimons, Megan and Athanassiadou, Maria and Kashem, Sheleem and Madge, Sam and others},
  journal = {Nature},
  volume  = {597},
  pages   = {672-677},
  year    = {2021},
  doi     = {10.1038/s41586-021-03854-z}
}

@article{Runge2019CausationEarth,
  title   = {Inferring causation from time series in Earth system sciences},
  author  = {Runge, Jakob and others},
  journal = {Nature Communications},
  volume  = {10},
  pages   = {2553},
  year    = {2019},
  doi     = {10.1038/s41467-019-10105-3}
}

@article{Schreiber2000,
  title   = {Measuring information transfer},
  author  = {Schreiber, Thomas},
  journal = {Physical Review Letters},
  volume  = {85},
  number  = {2},
  pages   = {461-464},
  year    = {2000},
  doi     = {10.1103/PhysRevLett.85.461}
}

@article{Sun2015OCSE,
  title   = {Optimal causation entropy for causal inference in complex systems},
  author  = {Sun, Jie and Taylor, Dane and Bollt, Erik M.},
  journal = {SIAM Journal on Applied Dynamical Systems},
  volume  = {14},
  number  = {1},
  pages   = {73-106},
  year    = {2015},
  doi     = {10.1137/140956166}
}

@article{almomani2020entropic,
  title   = {How entropic regression beats the outliers problem in nonlinear system identification},
  author  = {AlMomani, Abd AlRahman R. and Sun, Jie and Bollt, Erik M.},
  journal = {Chaos: An Interdisciplinary Journal of Nonlinear Science},
  volume  = {30},
  number  = {1},
  pages   = {013107},
  year    = {2020},
  doi     = {10.1063/1.5133386}
}

@article{fish2021entropic,
  title   = {Entropic regression with neurologically motivated applications},
  author  = {Fish, Jeremie and DeWitt, Alexander and AlMomani, Abd AlRahman R. and Laurienti, Paul J. and Bollt, Erik M.},
  journal = {Chaos: An Interdisciplinary Journal of Nonlinear Science},
  volume  = {31},
  number  = {11},
  pages   = {113105},
  year    = {2021},
  doi     = {10.1063/5.0039333}
}

@article{almomani2020erfit,
  title   = {ERFit: Entropic regression fit MATLAB package for data-driven system identification of underlying dynamic equations},
  author  = {AlMomani, Abd AlRahman and Bollt, Erik M.},
  journal = {arXiv preprint arXiv:2010.02411},
  year    = {2020},
  doi     = {10.48550/arXiv.2010.02411}
}

@article{achanta2012slic,
  title   = {SLIC Superpixels Compared to State-of-the-Art Superpixel Methods},
  author  = {Achanta, Radhakrishna and Shaji, Appu and Smith, Kevin and Lucchi, Aurelien and Fua, Pascal and Suesstrunk, Sabine},
  journal = {IEEE Transactions on Pattern Analysis and Machine Intelligence},
  volume  = {34},
  number  = {11},
  pages   = {2274-2282},
  year    = {2012},
  doi     = {10.1109/TPAMI.2012.120}
}

@article{diggans2022geometric,
  title   = {Geometric partition entropy: Coarse-graining a continuous state space},
  author  = {Diggans, Christopher Tyler and AlMomani, Abd AlRahman R.},
  journal = {Entropy},
  volume  = {24},
  number  = {10},
  pages   = {1432},
  year    = {2022},
  doi     = {10.3390/e24101432}
}

@article{diggans2025generalizing,
  title   = {Generalizing geometric partition entropy for the estimation of mutual information in the presence of informative outliers},
  author  = {Diggans, C. Tyler and AlMomani, Abd AlRahman R.},
  journal = {Chaos: An Interdisciplinary Journal of Nonlinear Science},
  volume  = {35},
  number  = {3},
  pages   = {033141},
  year    = {2025},
  doi     = {10.1063/5.0247397}
}

@article{kraskov2004estimating,
  title   = {Estimating mutual information},
  author  = {Kraskov, Alexander and St{\"o}gbauer, Harald and Grassberger, Peter},
  journal = {Physical Review E},
  volume  = {69},
  number  = {6},
  pages   = {066138},
  year    = {2004},
  doi     = {10.1103/PhysRevE.69.066138}
}

@article{frenzel2007partial,
  title   = {Partial mutual information for coupling analysis of multivariate time series},
  author  = {Frenzel, Stefan and Pompe, Bernd},
  journal = {Physical Review Letters},
  volume  = {99},
  number  = {20},
  pages   = {204101},
  year    = {2007},
  doi     = {10.1103/PhysRevLett.99.204101}
}

@article{sun2022data,
  title   = {Data-driven learning of Boolean networks and functions by optimal causation entropy principle},
  author  = {Sun, Jie and AlMomani, Abd AlRahman R. and Bollt, Erik M.},
  journal = {Patterns},
  volume  = {3},
  number  = {11},
  pages   = {100631},
  year    = {2022},
  doi     = {10.1016/j.patter.2022.100631}
}

@article{hochreiter1997lstm,
  title   = {Long Short-Term Memory},
  author  = {Hochreiter, Sepp and Schmidhuber, Juergen},
  journal = {Neural Computation},
  volume  = {9},
  number  = {8},
  pages   = {1735-1780},
  year    = {1997},
  doi     = {10.1162/neco.1997.9.8.1735}
}

@article{zhang2016mrms,
  title   = {Multi-Radar Multi-Sensor (MRMS) quantitative precipitation estimation: Initial operating capabilities},
  author  = {Zhang, Jian and Howard, Kenneth and Langston, Carrie and Kaney, Brian and Qi, Youcun and Tang, Lin and Grams, Heather and Wang, Yadong and Cocks, Stephen and Martinaitis, Steven and others},
  journal = {Bulletin of the American Meteorological Society},
  volume  = {97},
  number  = {4},
  pages   = {621-638},
  year    = {2016},
  doi     = {10.1175/BAMS-D-14-00174.1}
}

@book{Houze2014,
  author    = {Houze, Robert A.},
  title     = {Cloud Dynamics},
  edition   = {2},
  publisher = {Academic Press},
  address   = {Amsterdam},
  year      = {2014},
  doi       = {10.1016/C2009-0-63939-4}
}

@book{Markowski2010,
  author    = {Markowski, Paul and Richardson, Yvette},
  title     = {Mesoscale Meteorology in Midlatitudes},
  publisher = {Wiley-Blackwell},
  address   = {Chichester},
  year      = {2010},
  doi       = {10.1002/9780470682104}
}

@article{Adams1997Monsoon,
  author  = {Adams, David K. and Comrie, Andrew C.},
  title   = {The North American Monsoon},
  journal = {Bulletin of the American Meteorological Society},
  volume  = {78},
  number  = {10},
  pages   = {2197-2213},
  year    = {1997},
  doi     = {10.1175/1520-0477(1997)078<2197:TNAM>2.0.CO;2}
}

@book{Pearl2009,
  author    = {Pearl, Judea},
  title     = {Causality: Models, Reasoning, and Inference},
  edition   = {2},
  publisher = {Cambridge University Press},
  address   = {Cambridge},
  year      = {2009}
}

@book{Peters2017,
  author    = {Peters, Jonas and Janzing, Dominik and Sch{"o}lkopf, Bernhard},
  title     = {Elements of Causal Inference: Foundations and Learning Algorithms},
  publisher = {MIT Press},
  address   = {Cambridge, Massachusetts},
  year      = {2017}
}

@article{Murphy1993,
  author  = {Murphy, Allan H.},
  title   = {What Is a Good Forecast? An Essay on the Nature of Goodness in Weather Forecasting},
  journal = {Weather and Forecasting},
  volume  = {8},
  number  = {2},
  pages   = {281-293},
  year    = {1993},
  doi     = {10.1175/1520-0434(1993)008<0281:WIAGFA>2.0.CO;2}
}

@article{Gneiting2007,
  author  = {Gneiting, Tilmann and Raftery, Adrian E.},
  title   = {Strictly Proper Scoring Rules, Prediction, and Estimation},
  journal = {Journal of the American Statistical Association},
  volume  = {102},
  number  = {477},
  pages   = {359-378},
  year    = {2007},
  doi     = {10.1198/016214506000001437}
}

@article{Fawcett2006,
  author  = {Fawcett, Tom},
  title   = {An Introduction to ROC Analysis},
  journal = {Pattern Recognition Letters},
  volume  = {27},
  number  = {8},
  pages   = {861-874},
  year    = {2006},
  doi     = {10.1016/j.patrec.2005.10.010}
}

@article{Saito2015,
  author  = {Saito, Takaya and Rehmsmeier, Marc},
  title   = {The Precision-Recall Plot Is More Informative than the ROC Plot When Evaluating Binary Classifiers on Imbalanced Datasets},
  journal = {PLOS ONE},
  volume  = {10},
  number  = {3},
  pages   = {e0118432},
  year    = {2015},
  doi     = {10.1371/journal.pone.0118432}
}

@article{Gilleland2009,
  author  = {Gilleland, Eric and Ahijevych, David and Brown, Barbara G. and Casati, Barbara and Ebert, Elizabeth E.},
  title   = {Intercomparison of Spatial Forecast Verification Methods},
  journal = {Weather and Forecasting},
  volume  = {24},
  number  = {5},
  pages   = {1416-1430},
  year    = {2009},
  doi     = {10.1175/2009WAF2222269.1}
}

@book{Wilks2019,
  author    = {Wilks, Daniel S.},
  title     = {Statistical Methods in the Atmospheric Sciences},
  edition   = {4},
  publisher = {Elsevier},
  address   = {Amsterdam},
  year      = {2019}
}

@article{Jolliffe2016PCA,
  author  = {Jolliffe, Ian T. and Cadima, Jorge},
  title   = {Principal Component Analysis: A Review and Recent Developments},
  journal = {Philosophical Transactions of the Royal Society A},
  volume  = {374},
  number  = {2065},
  pages   = {20150202},
  year    = {2016},
  doi     = {10.1098/rsta.2015.0202}
}

@article{Brunton2016SINDy,
  author  = {Brunton, Steven L. and Proctor, Joshua L. and Kutz, J. Nathan},
  title   = {Discovering Governing Equations from Data by Sparse Identification of Nonlinear Dynamical Systems},
  journal = {Proceedings of the National Academy of Sciences},
  volume  = {113},
  number  = {15},
  pages   = {3932-3937},
  year    = {2016},
  doi     = {10.1073/pnas.1517384113}
}

@article{Barnett2009,
  author  = {Barnett, Lionel and Barrett, Adam B. and Seth, Anil K.},
  title   = {Granger Causality and Transfer Entropy Are Equivalent for Gaussian Variables},
  journal = {Physical Review Letters},
  volume  = {103},
  number  = {23},
  pages   = {238701},
  year    = {2009},
  doi     = {10.1103/PhysRevLett.103.238701}
}

@inproceedings{Runge2018CMI,
  author    = {Runge, Jakob},
  title     = {Conditional Independence Testing Based on a Nearest-Neighbor Estimator of Conditional Mutual Information},
  booktitle = {Proceedings of the Twenty-First International Conference on Artificial Intelligence and Statistics},
  series    = {Proceedings of Machine Learning Research},
  volume    = {84},
  pages     = {938-947},
  year      = {2018},
  publisher = {PMLR}
}

@article{Meinshausen2010Stability,
  author  = {Meinshausen, Nicolai and B{"u}hlmann, Peter},
  title   = {Stability Selection},
  journal = {Journal of the Royal Statistical Society: Series B},
  volume  = {72},
  number  = {4},
  pages   = {417-473},
  year    = {2010},
  doi     = {10.1111/j.1467-9868.2010.00740.x}
}

@article{Roberts2017CV,
  author  = {Roberts, David R. and Bahn, Volker and Ciuti, Simone and Boyce, Mark S. and Elith, Jane and Guillera-Arroita, Gurutzeta and Hauenstein, Severin and Lahoz-Monfort, Jos{\'e} J. and Schr{"o}der, Boris and Thuiller, Wilfried and Warton, David I. and Wintle, Brendan A. and Hartig, Florian and Dormann, Carsten F.},
  title   = {Cross-Validation Strategies for Data with Temporal, Spatial, Hierarchical, or Phylogenetic Structure},
  journal = {Ecography},
  volume  = {40},
  number  = {8},
  pages   = {913-929},
  year    = {2017},
  doi     = {10.1111/ecog.02881}
}

@inproceedings{Gao2017Mixture,
  author    = {Gao, Weihao and Kannan, Sreeram and Oh, Sewoong and Viswanath, Pramod},
  title     = {Estimating Mutual Information for Discrete-Continuous Mixtures},
  booktitle = {Advances in Neural Information Processing Systems},
  volume    = {30},
  year      = {2017}
}

@article{Mesner2021MixedCMI,
  author  = {Mesner, Octavio C. and Shalizi, Cosma R.},
  title   = {Conditional Mutual Information Estimation for Mixed, Discrete and Continuous Data},
  journal = {IEEE Transactions on Information Theory},
  volume  = {67},
  number  = {1},
  pages   = {464-484},
  year    = {2021},
  doi     = {10.1109/TIT.2020.3024886}
}

@article{Koster2004,
  author  = {Koster, Randal D. and Dirmeyer, Paul A. and Guo, Zhichang and Bonan, Gordon and Chan, Edmond and Cox, Peter and Gordon, C. T. and Kanae, Shinjiro and Kowalczyk, Eva and Lawrence, David and Liu, Ping and Lu, Chao-Han and Malyshev, Sergey and McAvaney, Bryant and Mitchell, Kenneth and Mocko, David and Oki, Taikan and Oleson, Keith and Pitman, Andy and Sud, Y. C. and Taylor, Christopher M. and Verseghy, Diana and Vasic, Ratko and Xue, Yongkang and Yamada, Tomohito},
  title   = {Regions of Strong Coupling between Soil Moisture and Precipitation},
  journal = {Science},
  volume  = {305},
  number  = {5687},
  pages   = {1138-1140},
  year    = {2004},
  doi     = {10.1126/science.1100217}
}

@article{Seneviratne2010,
  author  = {Seneviratne, Sonia I. and Corti, Thierry and Davin, Edouard L. and Hirschi, Martin and Jaeger, Eric B. and Lehner, Irene and Orlowsky, Boris and Teuling, Adriaan J.},
  title   = {Investigating Soil Moisture-Climate Interactions in a Changing Climate: A Review},
  journal = {Earth-Science Reviews},
  volume  = {99},
  number  = {3-4},
  pages   = {125-161},
  year    = {2010},
  doi     = {10.1016/j.earscirev.2010.02.004}
}

@article{Ford2015,
  author  = {Ford, Trent W. and Rapp, Anita D. and Quiring, Steven M. and Blake, Jennifer},
  title   = {Soil Moisture-Precipitation Coupling: Observations from the Oklahoma Mesonet and Underlying Physical Mechanisms},
  journal = {Hydrology and Earth System Sciences},
  volume  = {19},
  number  = {8},
  pages   = {3617-3631},
  year    = {2015},
  doi     = {10.5194/hess-19-3617-2015}
}

@article{almomani2020go,
  author  = {AlMomani, Abd AlRahman R. and Bollt, Erik M.},
  title   = {Go with the Flow, on Jupiter and Snow: Coherence from Model-Free Video Data without Trajectories},
  journal = {Journal of Nonlinear Science},
  volume  = {30},
  number  = {5},
  pages   = {2375-2404},
  year    = {2020},
  doi     = {10.1007/s00332-018-9470-1}
}

@article{almomani2021early,
  author  = {AlMomani, Abd AlRahman and Bollt, Erik},
  title   = {An Early Warning Sign of Critical Transition in the Antarctic Ice Sheet: A Data-Driven Tool for a Spatiotemporal Tipping Point},
  journal = {Nonlinear Processes in Geophysics},
  volume  = {28},
  number  = {1},
  pages   = {153-166},
  year    = {2021},
  doi     = {10.5194/npg-28-153-2021}
}

@article{diggans2023boltzmann,
  author  = {Diggans, C. Tyler and AlMomani, Abd AlRahman R.},
  title   = {Boltzmann-Shannon Interaction Entropy: A Normalized Measure for Continuous Variables with an Application as a Subsample Quality Metric},
  journal = {Chaos: An Interdisciplinary Journal of Nonlinear Science},
  volume  = {33},
  number  = {12},
  pages   = {123131},
  year    = {2023},
  doi     = {10.1063/5.0182349}
}

@article{Schultz2021,
  author  = {Schultz, Martin G. and Betancourt, Clara and Gong, Bing and Kleinert, Felix and Langguth, Michael and Leufen, Lukas H. and Mozaffari, Amir and Stadtler, Scarlett},
  title   = {Can Deep Learning Beat Numerical Weather Prediction?},
  journal = {Philosophical Transactions of the Royal Society A},
  volume  = {379},
  number  = {2194},
  pages   = {20200097},
  year    = {2021},
  doi     = {10.1098/rsta.2020.0097}
}

@misc{NOAAHRRRInventory,
  author       = {{NOAA National Centers for Environmental Prediction}},
  title        = {High-Resolution Rapid Refresh Product Inventory},
  howpublished = {National Centers for Environmental Prediction Central Operations},
  year         = {2026},
  url          = {https://www.nco.ncep.noaa.gov/pmb/products/hrrr/},
  note         = {Accessed 30 July 2026}
}

\end{document}